\documentclass[a4paper,11pt]{amsart}
\usepackage{geometry} 
\usepackage{amsmath,mathrsfs,amssymb,amsthm,bm,amsfonts}
\usepackage{color,graphicx,subfigure,epsfig}
\usepackage{times}
\allowdisplaybreaks

\newtheorem{theorem}{Theorem}[section]
\newtheorem{remark}{Remark}[section]

\newtheorem{lemma}{Lemma}[section]

\newtheorem{assumption}{Assumption}[section]
\newtheorem{alg}{Algorithm}
\newcommand{\p}{\partial}
\newcommand{\f}{\frac}

\begin{document}
\title[Hybrid mixed discontinuous Galerkin method]{Hybrid mixed discontinuous Galerkin finite element method for incompressible wormhole propagation problem}
\thanks{Corresponding author: Jiansong Zhang.
\\Zhang's work was supported by the Natural Science Foundation of Shandong Province (ZR2019MA015),  and the Fundamental Research Funds for the Central Universities  (22CX03020A). Zhu’s work was partially supported by the National Council for Scientific and Technological Development of Brazil (CNPq).}

\author{Jiansong Zhang, Yun Yu, Jiang Zhu, Yue Yu and Rong Qin}

\begin{abstract}
Wormhole propagation plays a very important role in the product enhancement of oil and gas reservoir. A new combined hybrid mixed finite element method is proposed to solve incompressible wormhole propagation problem with discontinuous Galerkin finite element procedure, in which, the new hybrid mixed finite element algorithm is established for pressure equation, while the discontinuous Galerkin finite element method is considered for concentration equation, and then the porosity function is computed straightly by the approximate value of the concentration. This new combined method can keep local mass balance, meantime it also keeps the boundedness of the porosity. The convergence of the proposed method is analyzed and the optimal error estimate is derived. Finally, numerical examples are presented to verify the validity of the algorithm and the correctness of the theoretical results.
\end{abstract}



\date{} 
\address{Jiansong Zhang: College of Science,  China University of Petroleum, Qingdao 266580, China}
\email{jszhang@upc.edu.cn}

\address{Yun Yu: College of Science, China University of Petroleum, Qingdao 266580, China}
\email{yuyun19970321@163.com}

\address{Jiang Zhu: Laborat\'orio Nacional de Computa\c{c}\~ao Cient\'{\i}fica, MCTI\\
        Avenida Get\'ulio Vargas 333, 25651-075 Petr\'opolis, RJ, Brazil}
\email{jiang@lncc.br}

\address{Yue Yu: College of Science, China University of Petroleum, Qingdao 266580, China}
\email{m18766215811@163.com}

\address{Rong Qin: College of Science, China University of Petroleum, Qingdao 266580, China}
\email{qr-920diana@126.com}

\maketitle
\section{Introduction}
The acid treatment of carbonate reservoirs is a widely practiced oil and gas well stimulation technique. In fact, when acids are injected into oil production wells, chemical reactions cause the dissolution of the material near the wellbore to result into flow channels. Such flow channels look like worm holes that they are usually called wormholes.  Because of its important role in the product enhancement of oil and gas reservoir,  the wormhole propagations have been a topic of key interest for research during recent decades. The theoretical researches on numerical  methods for these problems have extensive practicability and important significance. 

Here, we will construct a new combined numerical procedure to solve the incompressible wormhole propagation problem which is usually described by the following nonlinear partial differential equations  (see \cite{zhao1,rui2,sun1}):
\begin{equation}\label{eq1}
\left\{\begin{aligned}
&\frac{\partial\phi}{\partial t}=\frac{\alpha k_ca_v(c_f-c_s)}{\rho_s},\quad x\in \Omega,\quad 0 \le t \le T,\\
&\frac{\partial\phi}{\partial t}+\nabla\cdot \mathbf{u} =f,\quad x\in \Omega,\quad 0 \le t \le T, \\
&\mathbf{u}=\frac{-k(\phi)}{\mu}\nabla p,\quad x\in \Omega,\quad 0 \le t \le T, \\
&\frac{\partial(\phi c_f)}{\partial t}+\nabla\cdot(\mathbf{u}c_f)=\nabla\cdot(\phi \mathbf{D}(\mathbf{u})\nabla c_f)+k_ca_v(c_s-c_f)+f_Ic_I-f_Pc_f,
\end{aligned}\right.
\end{equation}
and the  corresponding initial-boundary conditions are considered as follows:
\begin{equation}\label{boundary}
\left\{\begin{aligned}
&\phi(x,0)=\phi_{0}(x),\quad c_{f}(x,0)=c^{0}_{f}(x),\quad x\in \Omega,\\
&\mathbf{u}\cdot \mathbf{n}=0,\quad (\phi \mathbf{D}(\mathbf{u})\nabla c_f-c_f\mathbf{u})\cdot \mathbf{n}=0,\quad x\in \partial\Omega,\quad 0 \le t \le T,
\end{aligned}\right.
\end{equation}
where $\Omega\subset R^d(d=2,3)$ denotes a bounded polygonal/polyhedral domain;  $\alpha$ is the dissolving constant of the acid; $\rho_s$ is the density of the rock; $\mathbf{n}$ is the unit outward normal vector  to $\partial\Omega$; the functions $p$ and $\mathbf{u}$ denote the pressure and Darcy velocity; $\phi$ and $k$ are the porosity and permeability of rocks, and $\mu$ is the viscosity of fluid;  $a_v$ is the interfacial area
available for reaction; $f$, $f_I$ and $f_P$ are the external volumetric flow rate, the injection flow rate and the production flow rate, respectively; $c_f$, $c_s$ and $c_I$ are the concentrations of acid in the fluid phase, the fluid-solid interface and the injected flow, respectively; Diffusion coefficient $\mathbf{D}(\mathbf{u})=\phi[d_mI+|\mathbf{u}|(d_IE(\mathbf{u})+d_tE^{\perp}(\mathbf{u}))]$ comes from two aspects: small molecule diffusion of oil field scale problem, and speed-related diffusion in petroleum engineering, where the matrix $E(\mathbf{u})=(u_iu_j/|\mathbf{u}|^2)_{d \times d}$ and $E^{\perp}(\mathbf{u}) = I-E(\mathbf{u})$;  $k_c$ is the local mass-transfer coefficient.
In the case of first order kinetic reaction, the concentrations $c_s$ and $c_f$ satisfy the relationship:
\begin{equation}\label{5}
c_s=\frac{c_f}{1+k_s/k_c},
\end{equation}
where $k_s$ is the kinetic constant for reaction. The quantitative relationship between the coefficients $k$, $v$ and $\phi$ is as follows:
\begin{equation}\label{6}
\begin{aligned}
&\frac{k}{k_0}=\frac{\phi}{\phi_0}(\frac{\phi(1-\phi_0)}{\phi_0(1-\phi)})^2
,\quad \frac{a_v}{a_0}=\frac{1-\phi}{1-\phi_0},
\end{aligned}
\end{equation}
where $k_0$, $a_0$, and $\phi_0$ are the initial values for $k$, $a_v$, $\phi$.

Many numerical methods have been constructed for simulating the wormhole propagation.  In  \cite{sun1}, Kou etc. used a classical mixed finite element procedure to establish  a  fully conservative method for incompressible wormhole problem;  And then, they considered  a parallel algorithm for wormhole problem in \cite{sun2} under the  Darcy-Brinkman-Forchheimer framework.  In  \cite{rui2,rui1}, Rui and Li studied the block-centered finite difference methods with or without the method of characteristics for wormhole propagation.  Zhang etc. established a combined splitting mixed finite element method for  compressible wormhole propagation with the method of the characteristics in \cite{zhang0}.  And Guo etc. studied the local discontinuous Galerkin finite element method for incompressible problem in \cite{guo1}.

Generally speaking, the mixed finite element methods can be used to obtain more  accurate approximation of the velocity function. However, the classical mixed element methods  usually result  into some saddle point problems, in which the coefficient matrix of the mixed system loses the symmetric positive definite property and the finite element spaces require the LBB condition. By introducing the Lagrange multiplier, a symmetric and positive definite system is obtained by hybrid mixed element method  ( \cite{ref5,ref6,ref7,ref8,zhang,zhu5} ).  Therefore, the hybrid technique will be considered for the pressure and velocity equations.  In the procedure, the velocity and pressure are eliminated by use of the Lagrange multiplier,  so that they can be solved  element by element. In addition, the resulted global mixed system only involves the degrees of freedom with the Lagrange multiplier, so this technique can significantly improve the computational efficiency.

As we know, the concentration equation is usually characterized as convection-dominant, the traditional Galerkin finite element method is not well applied, in particular for the discontinuous problems. The discontinuous Galerkin (DG) methods  in \cite{ref9,ref10,ref11} were introduced. They have several advantages over other types of finite element methods. For example, test functions across the finite element
interfaces do not explicitly impose continuity constraints. As a result, the finite element spaces allow highly nonuniform and unstructured meshes. These methods have been more and more widely used, such as \cite{ref12,ref13}. Zhu and his coauthors studied the discontinuous Galerkin finite element methods for  nonlinear reaction diffusion equations in \cite{zhu3,zhu4}, and then extended these techniques to the dipolar Bose-Einstein condensation  and Lagrangian compressible Euler equations in \cite{zhu1,zhu2}.  Based on the traditional mixed finite element methods and DGFE methods, the combined mixed DG methods were proposed to solve the compressible and incompressible miscible displacement problems in \cite{ref14,ref15,ref16,ref20}.  However,  there is little research on the discontinuous Galerkin methods for wormhole propagation.

The focus of this article is to combine the discontinuous Galerkin finite element method with the hybrid mixed technique  to simulate the incompressible wormhole propagation. In the combined method, a new hybrid mixed finite element (HMFE) procedure is constructed to solve pressure equation,  and the symmetric interior penalty discontinuous Galerkin (SIPDG) procedure is proposed to solve the concentration equation,  then  the porosity  is computed straightly by the approximated value of the concentration. Compared with other existing combined methods, the proposed method not only keeps mass balance locally, but also keeps the boundedness of the porosity, especially it can deal well with the discontinuous case. The consistency and stability of the proposed method are analyzed, and then the corresponding error estimate is given  under the case that the diffusion coefficient includes the molecular diffusion and dispersion, unlike the ones in \cite{rui2,rui1,guo1} where only  molecular diffusion was considered. Finally, numerical examples are presented to verify the validity of the algorithm and the correctness of the theoretical results.

For the convenience  of analysis,  we  make the following assumptions:
\begin{assumption}\label{a1}
Assume that the parameters
$\mu$, $k_c$, $k_s$, $\alpha$, $\rho_s$ are positive constants, and that $\phi_0$, $\frac{k(\phi)}{\mu}$ and $f(\cdot,t)$ are bounded as follows:
\begin{equation}
0<a_{\ast}\leq\frac{k(\phi)}{\mu}\leq a^{\ast},\quad 0<\phi_{0}< 1,\quad
| f(\cdot,t)|\leq C, 
\end{equation}
where  $a_{\ast}$, $a^{\ast}$ and $C$ are some positive constants. And we also assume that the diffusion coefficient $\mathbf{D}(\mathbf{u})$ satisfies the uniformly positive definiteness and Lipschitz continuousness 
\begin{equation}
\begin{aligned}\nonumber
&\mathbf{D}(\mathbf{u})\nabla c\cdot\nabla c\geq D_{*}| \nabla c|^2
\end{aligned}
\end{equation}
and 
\begin{equation}\label{2222}
\begin{aligned}
&\| \mathbf{D}(\mathbf{u})-\mathbf{D}(\mathbf{v})\|_{[L^2]^d}\leq D^{*}\| \mathbf{u}-\mathbf{v}\|_{[L^2]^d},
\end{aligned}
\end{equation}
where $D*$ and $D_{*}$ are two positive constants independent of $\mathbf{u}$ and $\mathbf{v}$ and $c$.
\end{assumption}

\begin{assumption}\label{a2}
Assume that the solution $(p, \mathbf{u}, c_f)$ of the system \eqref{eq1} has the regularities as follows:
\begin{equation}
\begin{aligned}
&({a})\quad p\in L^2(0,T;H^{k+2}(\Omega)),\quad \phi_{0}\leq \phi\leq C_{1}<1,\\
&({b})\quad \mathbf{u}\in L^\infty(0,T;H^{k+1}(\Omega))\cap L^\infty(0,T;L^\infty(\Omega)),\\
&({c})\quad c_f\in H^{1}(0,T;H^{k+1}(\Omega))\cap L^{\infty}(0,T;W^{1,\infty}(\Omega)),\\
&({d})\quad\frac{\partial c_f}{\partial t}\in L^2(0,T;H^{k+1}(\Omega))\cap L^{\infty}(0,T;L^{\infty}(\Omega)).
\end{aligned}
\end{equation}
\end{assumption}

Moreover,  we only consider the homogeneous boundary
condition case in this article for simplification. For nonhomogeneous boundary value problem, we can use some simple technique to transform it into homogeneous boundary value problem, so our method proposed later is still valid.

\section{The formulation of HMDG method}

In order to illustrate our method, we first give a uniform regular partition of $\Omega$, denoted by $\mathcal{T}_h = \{K_1, K_2, . . . , K_n \}$. We denote $\mathcal {E}_h=\cup_{K\in\mathcal{T}_h}\{e|e\in\partial K\}$ to be the set of all cell edges and $\mathcal{E}_h^i=\mathcal{E}_h\backslash \partial\Omega$ as all the interior ones. Furthermore, let $h_e = diam(e)$ for all $e\in\mathcal{E}_h$.
Introduce the following piecewise Sobolev spaces associated with ${\mathcal{T}_h}$
\begin{equation}\nonumber
H^s(\mathcal{T}_h)=\{v\in L^2(\Omega):v|_K\in H^s(K),K\in\mathcal{T}_h\},\quad s\ge 0.
\end{equation}

We also define the following spaces on $\mathcal{E}_h$:
\begin{equation}\nonumber
\begin{aligned}
&L^2(\mathcal{E}_h)=\{v\in L^2(e),\forall e\in\mathcal{E}_h\},
&L^2(\mathcal{E}_h^i)=\{v\in L^2(e),\forall e\in\mathcal{E}_h^i\}.
\end{aligned}
\end{equation}

For $e\in\mathcal{E}_h^i$, denote by $n_e$ a fixed unit normal direction. For $e\in\partial\Omega$,  $n_e=\mathbf{n}$. We define averages $\{\cdot\}$ and jumps $[\cdot]$:
\begin{equation}\nonumber
\{v\}=\frac{1}{2}[(v|_{K^{i}})|_e+(v|_{K^{j}})|_e],\quad[v]=(v|_{K^{i}})|_e-(v|_{K^{j}})|_e,\quad \textrm{on}\quad e\in\mathcal{E}_h^i.
\end{equation}
In particular, if $e\in\partial\Omega$, $\{v\}=v|_e,\quad[v]=v|_e$. Meantime, we define inner products as follows:
\begin{equation}\nonumber
(\cdot,\cdot)_{\mathcal{T}_h}=\sum\limits_{K\in\mathcal{T}_h}(\cdot,\cdot)_K,\quad \langle\cdot,\cdot\rangle_{\mathcal{E}_h}=\sum\limits_{e\in\mathcal{E}_h}\langle\cdot,\cdot\rangle_e,\quad \langle\cdot,\cdot\rangle_{\mathcal{E}_h^i}=\sum\limits_{e\in\mathcal{E}_h^i}\langle\cdot,\cdot\rangle_e,
\end{equation}
and the norms $\|\cdot\|_{\mathcal{T}_{h}}=\sqrt{(\cdot,\cdot)_{\mathcal{T}_h}}$ and  $|\cdot|_{\p\mathcal{T}_{h}}=\sqrt{\langle\cdot,\cdot\rangle_{\mathcal{E}_h}}$.

Introduce the discrete approximate spaces denoted by $\Psi_h$,$\Lambda_h$,$\Pi_h$ and $\Sigma_h$ as follows:
\begin{equation}\nonumber
\begin{aligned}
&\Psi_h=\{v\in H^k(\mathcal{T}_h):v|_K\in P_k(K),K\in\mathcal{T}_h\},\\
&\Lambda_h=\{v\in L^2(\mathcal{T}_h):v|_K\in P_k(K),K\in\mathcal{T}_h\},\\
&\Pi_h=\{v\in [H^k(\mathcal{T}_h)]^d:v|_K\in RT_k(K),K\in\mathcal{T}_h\},\\
&\Sigma_h=\{v\in L^2(\mathcal{E}_h):v|_e\in P_k(e),e\in\mathcal{E}_h\},
\end{aligned}
\end{equation}
where $P_k(K)$,$P_k(e)$ are the spaces of polynomial functions of degree at most $k$ for each $K\in\mathcal{T}_h$ and each $e\in\mathcal{E}_h$, respectively, $RT_k(K)=[P_k(K)]^d\oplus xP_k(K)$  denotes the Raviart-Thomas mixed finite element space as in \cite{ref6,ref4,ref17}.

Set
\[
\kappa=\frac{k_ck_sa_0}{(k_c+k_s)(1-\phi_0)}.
\]
Using \eqref{5} and \eqref{6}, we can rewrite \eqref{eq1} in the following equivalent form:
\begin{equation}\label{model:2}
\begin{aligned}
&\frac{\partial\phi}{\partial t}=\frac{\alpha \kappa}{\rho_s}(1-\phi)c_f,\\
&\frac{\alpha \kappa}{\rho_s}(1-\phi)c_f+\nabla\cdot\mathbf{u} =f,\quad \mathbf{u}=\frac{-k(\phi)}{\mu}\nabla p,\\
&\frac{\partial(\phi c_f)}{\partial t}+\nabla\cdot(\mathbf{u}c_f)=\nabla\cdot(\phi \mathbf{D}(\mathbf{u})\nabla c_f)- \kappa(1-\phi)c_f+f_Ic_I-f_Pc_f,.
\end{aligned}
\end{equation}

Next, we will formulate our method for wormhole propagation.

For the discretization of the porosity, we consider the similar technique as in \cite{sun1}. The discrete porosity is point-wise defined and can be stated as follows:
\begin{equation}\label{semi}
\begin{aligned}
&\frac{\partial\phi_h}{\partial t}=\frac{\alpha \kappa}{\rho_s}(1-\phi_h)\bar{c}_h,
\end{aligned}
\end{equation}
where $\bar{c}_h=\max(0,\min(c_h,1))$, $c_{h}$ is a given approximation of the concentration $c_{f}$.

\subsection{Hybrid mixed finite element scheme for pressure and velocity}
In this subsection, we give the hybrid mixed finite element (HMFE) method for pressure and velocity. which can be written as below:
\begin{alg}
For given approximate values of $\phi_{h}$ and $c_{h}$, seek $(p_h, \mathbf{u}_h, \lambda_h) \in  \Lambda_h \times \Pi_h \times \Sigma_h$ such that
\begin{equation}\label{hmfm}
\begin{aligned}
&\sum\limits_{K\in\mathcal{T}_h}(\frac{\alpha \kappa}{\rho_s}(1-\phi_h)\bar{c}_h,v_h)_K+\sum\limits_{K\in\mathcal{T}_h}(\nabla\cdot \mathbf{u}_h,v_h)_K=\sum\limits_{K\in\mathcal{T}_h}(f,v_h)_K,\quad\forall v_h \in \Lambda_h,\\
&\sum\limits_{K\in\mathcal{T}_h}(a(\phi_h)\mathbf{u}_h,\omega_h)_K-\sum\limits_{K\in\mathcal{T}_h}(p_h,\nabla\cdot\omega_h)_K+\sum\limits_{e\in\mathcal{E}_h}\langle\lambda_h,[\omega_h]\cdot n_e\rangle_e=0,\quad\forall \omega_h \in \Pi_h,\\
&\sum\limits_{e\in\mathcal{E}_h^i}\langle\mu_h,[\mathbf{u}_h]\cdot n_e\rangle_e=0,\quad\forall \mu_h \in \Sigma_h,
\end{aligned}
\end{equation}
where $a(\phi_h)=\frac{\mu}{k(\phi_h)}$ and $\bar{c}_h=\max(0,\min(c_h,1))$.
\end{alg}

Define the bilinear form:
\[
\begin{aligned}
B_{\mathbf{u}}((\mathbf{u}_h,p_h,\lambda_h),(\omega_h,v_h,\mu_h))=:&(\mathbf{u}_h,\nabla v_h)_{\mathcal{T}_h}+(\nabla p_h,\omega_h)_{\mathcal{T}_h}+(a(\phi_h)\mathbf{u}_h,\omega_h)_{\mathcal{T}_h}\\
&+\langle\lambda_h-p_h,[\omega_h]\cdot n_e\rangle_{\mathcal{E}_h}+\langle[\mathbf{u}_h]\cdot n_e,\mu_h-v_h\rangle_{\mathcal{E}_h^i},\\
B_{\phi}(c_h;\phi_h,v_h):=&-(\frac{\alpha \kappa}{\rho_s}(1-\phi_h)\bar{c}_h,v_h)_{\mathcal{T}_h}.
\end{aligned}
\]
We can rewrite \eqref{hmfm} into the following equivalent form:
\begin{alg}[HMFE Algorithm]  For given $\phi_h$ and $c_h$, find $(\mathbf{u}_h, p_h, \lambda_h) \in  \Pi_h \times\Lambda_h \times \Sigma_h$ such that
\begin{equation}
\begin{array}{c}
B_{\phi}(c_h;\phi_h,v_h)+B_{\mathbf{u}}((\mathbf{u}_h,p_h,\lambda_h),(\omega_h,v_h,\mu_h))=-(f,v_h)_{\mathcal{T}_h}, \\[.1in]
 \forall ( \omega_h, v_h,\mu_h) \in \Pi_h \times \Lambda_h \times \Sigma_h.
 \end{array}
\end{equation}
\end{alg}
\subsection{DGFE method for the concentration}
Due to the flexibility of the discontinuous Galerkin finite element method in constructing feasible local-shape function spaces and the advantage in capturing non-smooth or oscillatory solutions effectively, we consider it to be applied for the concentration.

Define the bilinear form
\[
\begin{aligned}
B_c(c_h,z_h):=&\sum\limits_{K\in\mathcal{T}_h}\int_{K}(\phi_h\mathbf{D}(\mathbf{u}_h)\nabla c_h-\mathbf{u}_hc_h)\nabla z_hdx-\sum\limits_{e\in\mathcal{E}_h^i}\int_{e}\{(\phi_h\mathbf{D}(\mathbf{u}_h)\nabla z_h-\mathbf{u}_hz_h)\cdot n_e\}[c_h]ds\\
&-\sum\limits_{e\in\mathcal{E}_h^i}\int_{e}\{(\phi_h\mathbf{D}(\mathbf{u}_h)\nabla c_h-\mathbf{u}_hc_h)\cdot n_e\}[z_h]ds+J_0^\gamma(c_h,z_h),
\end{aligned}
\]
where $J_0^\gamma(c_h,z_h)$ denotes  the penalty term defined by
\[
J_0^\gamma(c_h,z_h)=\sum\limits_{e\in\mathcal{E}_h^i}\int_{e}\frac{\gamma}{h_e^\beta}[c_h][z_h],
\]
where $\gamma$ is called penalty parameter and bounded below by a large enough constant, and $\beta$ denotes some positive constant.

Now we can reach the SIPDG method for the concentration equation.
\begin{alg}[SIPDG Algorithm]  For given $\mathbf{u}_h$ and $\phi_h$, find $c_h\in\Psi_h$ such that
\begin{equation}
(\frac{\partial(\phi_hc_h)}{\partial t},z_h)_{\mathcal{T}_h}+B_c(c_h,z_h)+(\kappa(1-\phi_h)c_{h},z_h)_{\mathcal{T}_h}=(f_Ic_I-f_Pc_h,z_h)_{\mathcal{T}_h},\quad \forall  z_h\in\Psi_h.
\end{equation}
\end{alg}

\subsection{The combined HMDG method}
Here, we will present the new combined  SIPDG finite element method for incompressible wormhole problem with hybrid mixed finite element procedure.
\begin{alg}[HMDG Algorithm]\label{alg4}
For the given initial value functions $\phi_{0}$ and $c^{0}_{f}$, find $(c_h, p_h, \mathbf{u}_h, \lambda_h) \in \Psi_h \times \Lambda_h \times \Pi_h \times \Sigma_h$ such that
\begin{equation}\label{101}
\begin{aligned}
&(\textrm{a})\quad\frac{\partial\phi_h}{\partial t}=\frac{\alpha \kappa}{\rho_s}(1-\phi_h)\bar{c}_h,\quad \bar{c}_h=\max(0,\min(c_h,1)),\\
&(\textrm{b})\quad B_{\phi}(c_h;\phi_h,v_h)+B_{\mathbf{u}}((\mathbf{u}_h,p_h,\lambda_h),(\omega_h,v_h,\mu_h))=-(f,v_h)_{\mathcal{T}_h},\\
&(\textrm{c})\quad(\frac{\partial(\phi_hc_h)}{\partial t},z_h)_{\mathcal{T}_h}+B_c(c_h,z_h)+(\kappa(1-\phi_h)c_{h},z_h)_{\mathcal{T}_h}=(f_Ic_I-f_Pc_h,z_h)_{\mathcal{T}_h}.
\end{aligned}
\end{equation}
\end{alg}
\begin{theorem}
Algorithm \ref{alg4} is consistent. That is, define $\lambda = p$,  then the solution $(p, u, \lambda,\phi, c)$ of system \eqref{model:2} satisfies  \eqref{101}.
Conversely, if $(p_h, u_h, \lambda_h,\phi_h ,c_h)$ satisfies \eqref{101}, then $(p_h, u_h, \lambda_h,\phi_h, c_h)$ are the solutions of problem \eqref{model:2}.
\end{theorem}
\begin{proof}
 Substituting the weak solution $(p, u, \phi, c)$ of problem \eqref{model:2} into the first equation of \eqref{101} with $\omega_h=\mu_h=0$, we can get
 \[
\begin{aligned}
-(\frac{\alpha \kappa}{\rho_s}(1-\phi)c_f,v_h)_{\mathcal{T}_h}-(\nabla\cdot\mathbf{u},v_h)_{\mathcal{T}_h}=-(f,v_h)_{\mathcal{T}_h}.
\end{aligned}
\]

Next, we test with $v_h=\omega_h=0$ and get the equation
\[
\langle\mathbf{u}\cdot n_e,\mu_h\rangle_{\mathcal{E}_h^i}=0,
\]
so the normal component of the flux $\mathbf{u}$ is continuous at element
interfaces.

Now we will prove the consistency between the model problem \eqref{model:2} and HMDG Algorithm \eqref{101}.

Firstly, we prove that the solution of \eqref{model:2} also solves \eqref{101}. For this, let $z$ be an element in $\Psi_h$. We multiply the third equation of \eqref{model:2} by $z$ and integrate
on one element $K$:

\begin{equation}\nonumber
\begin{aligned}
&\int_{K}\frac{\partial(\phi c_f)}{\partial t}zdx+\int_{K}(\phi\mathbf{D}(\mathbf{u})\nabla c_f-\mathbf{u}c_f)\nabla zdx-\int_{\partial K}(\phi\mathbf{D}(\mathbf{u})\nabla c_f-\mathbf{u}c_f)\cdot n_Kzds\\
=&\int_{K}(f_Ic_I-f_Pc_f-\kappa(1-\phi)c_{f})z)dx.
\end{aligned}
\end{equation}

Summing it over all elements and using \eqref{boundary}, we observe that
\[
\begin{aligned}
\sum\limits_{K\in\mathcal{T}_h}\int_{\partial K}(\phi\mathbf{D}(\mathbf{u})\nabla c_f-\mathbf{u}c_f)\cdot n_ezds
=\sum\limits_{e\in\mathcal{E}_h^i}\int_{e}(\phi\mathbf{D}(\mathbf{u})\nabla c_f-\mathbf{u}c_f)\cdot n_e[z]ds.
\end{aligned}
\]

By the regularities of the solution $\phi$, $u$ and $c$, we have
\begin{equation}\nonumber
\begin{aligned}
&(a)\quad(\phi\mathbf{D}(\mathbf{u})\nabla c_f-\mathbf{u}c_f)\cdot n_e[z]=\{(\phi\mathbf{D}(\mathbf{u})\nabla c_f-\mathbf{u}c_f)\cdot n_e\}[z],\\
&(b)\quad[c_f]=0.
\end{aligned}
\end{equation}
Therefore, we obtain the second equation of the scheme \eqref{101}.

Conversely, take $z\in H^{1}(\Omega)$ and $c_f\in H^{1}(\Omega)\cap\Psi_h$. Then \eqref{101} reduces to
\begin{equation}\nonumber
\begin{aligned}
\sum\limits_{K\in\mathcal{T}_h}\int_{K}\frac{\partial(\phi c_f)}{\partial t}z_hdx+\sum\limits_{K\in\mathcal{T}_h}\int_{ K}(\phi\mathbf{D}(\mathbf{u})\nabla c_f-\mathbf{u}c_f)zdx
=\sum\limits_{K\in\mathcal{T}_h}\int_{K}(f_Ic_I-f_Pc_f-\kappa(1-\phi)c_{f})z)dx.
\end{aligned}
\end{equation}

So for all $K\in\mathcal{T}_h$, we can obtain
\begin{equation}\nonumber
\begin{aligned}
\frac{\partial(\phi c_f)}{\partial t}+\nabla\cdot(\mathbf{u}c_f)=\nabla\cdot(\phi\mathbf{D}(\mathbf{u})\nabla c_f)+f_Ic_I-f_Pc_f-\kappa(1-\phi)c_{f},\quad\quad \textrm{in}\quad K.
\end{aligned}
\end{equation}

Finally, let $K_1$ and $K_2$ to be two adjacent elements,  and $e=\p K_1\cap\p K_2$. Take $z \in C^{\infty}_0 (K_1 \cup K_2 )$ and extend it by zero over the rest of the domain. Integrating by parts in the second equation of \eqref{model:2}, we can get
\begin{equation}\nonumber
\begin{aligned}
&(\frac{\partial(\phi c_f)}{\partial t},z)_{K_1 \cup K_2}+(\phi \mathbf{D}(\mathbf{u})\nabla c_f-\mathbf{u}c_f,\nabla z)_{K_1 \cup K_2}-\langle[(\phi \mathbf{D}(\mathbf{u})\nabla c_f-\mathbf{u}c_f)\cdot n_e],z\rangle_e\\
=&(f_Ic_I-f_Pc_f-\kappa(1-\phi)c_{f},z)_{K_1 \cup K_2}.
\end{aligned}
\end{equation}

On the other hand, \eqref{101} reduces to
\begin{equation}\nonumber
\begin{aligned}
(\frac{\partial(\phi c_f)}{\partial t},z)_{K_1 \cup K_2}+(\phi \mathbf{D}(\mathbf{u})\nabla c_f-\mathbf{u}c_f,\nabla z)_{K_1 \cup K_2}=(f_Ic_I-f_Pc_f-\kappa(1-\phi)c_{f},z)_{K_1 \cup K_2}.
\end{aligned}
\end{equation}
Hence, we have
\[
\langle[(\phi D(\mathbf{u})\nabla c_f-\mathbf{u}c_f)\cdot n_e],z\rangle_e=0,\quad\forall z \in C^{\infty}_0 (K_1 \cup K_2 ).
\]
Since this holds for all $e$, it implies that
$(\phi \mathbf{D}(\mathbf{u})\nabla c_f-\mathbf{u}c_f)\cdot n_e=0$ on $\partial\Omega$ and $\nabla\cdot(\phi \mathbf{D}(\mathbf{u})\nabla c_f-\mathbf{u}c_f)\in L^2(\Omega)$, hence we have
\begin{equation}
\begin{aligned}
&\frac{\partial(\phi c_f)}{\partial t}+\nabla\cdot(\mathbf{u}c_f-\phi \mathbf{D}(\mathbf{u})\nabla c_f)=f_Ic_I-f_Pc_f-\kappa(1-\phi)c_{f}.
\end{aligned}
\end{equation}
\end{proof}

We can easily show that the discrete solution $\phi_{h}$ of $\phi$ satisfies the following boundedness.
\begin{theorem}[The boundedness of porosity]\label{phih}
For any time $t\in (0,T]$,
the approximate porosity $\phi_h$ is bounded, that is,
\begin{equation}\label{phihbou1}
\phi_0\leq\phi_h\leq1-(1-\phi_0)e^{-\eta t}<1,
\end{equation}
\begin{equation}\label{phihbou2}
\begin{aligned}
&0\leq\frac{\partial\phi_h}{\partial t}\leq\frac{\alpha \kappa}{\rho_s},
\end{aligned}
\end{equation}
where $\eta=\frac{\alpha \kappa}{\rho_s}$ and $\phi_0>0$.
\end{theorem}
\begin{proof}
First, we can rewrite \eqref{semi} as the following integral form
\begin{equation}\nonumber
\begin{aligned}
&\int_{\phi_0}^{\phi_h}\frac{1}{1-\phi_h}d\phi=\int_{0}^{t}\eta\bar{c}_hd\tau,\\
\end{aligned}
\end{equation}
where $\eta>0$.
We deduce that
\begin{equation}\nonumber
\begin{aligned}
&\ln\frac{1-\phi_h}{1-\phi_0}=-\eta\bar{c}_ht.
\end{aligned}
\end{equation}
Notice that
\begin{equation}\nonumber
\begin{aligned}
\phi_h=1-(1-\phi_0)e^{-\eta\bar{c}_ht}\leq1-(1-\phi_0)e^{-\eta t}.
\end{aligned}
\end{equation}
It is easily seen that the approximate value of the porosity increases with $t$, and $\phi_h=\phi_0$ at $t = 0$, so we can get the estimate \eqref{phihbou1}. The estimate \eqref{phihbou2} is reached by \eqref{phihbou1} and \eqref{semi}.
\end{proof}

For HMDG Algorithm, we have the following main convergence theorem.
\begin{theorem}\label{thm1}
Under the assumptions \eqref{a1} and \eqref{a2}, for $t>0$, there is a priori  error estimate as follows:
\begin{equation}\label{3191}
\left\{\begin{split}
&(\textrm{a})\quad\| \phi-\phi_h\|_{\mathcal{T}_h}+\| \mathbf{u}-\mathbf{u}_h\|_{\mathcal{T}_h}+\| c_f-c_h\|_{\mathcal{T}_h}
\leq Ch^s(\| c_f\|_{H^1(0,T;H^s(\mathcal{T}_h))}+\| p\|_{L^2(0,T;H^{s+1}(\mathcal{T}_h))}),\\
&(\textrm{b})\quad\| \nabla(p_h-\Pi_h p)\|_{\mathcal{T}_h}+h^{-\frac{1}{2}}\| \lambda_h-p\|_{\mathcal{T}_h}
\leq Ch^s(\| c_f\|_{H^1(0,T;H^s(\mathcal{T}_h))}+\| p\|_{L^2(0,T;H^{s+1}(\mathcal{T}_h))}),\\
&(\textrm{c})\quad\| p-p_h\|_{\mathcal{T}_h}\leq Ch^s(\| c_f\|_{H^1(0,T;H^s(\mathcal{T}_h))}+\| p\|_{L^2(0,T;H^{s+1}(\mathcal{T}_h))}),
\end{split}\right.
\end{equation}
where when $d=2$, $1\leq s\leq k+1$; when $d=3$, $3/2\leq s\leq k+1$.
\end{theorem}

\section{Some important projections and lemmas}
In this section, we will give some important projection operators and approximate properties, which is used to show the convergence theorem of our proposed method.

Firstly, we introduce the following norms with respect to the bilinear form $B_{\mathbf{u}}$:
\[
\|(\omega,v,\mu)\||^2_B:=\|\omega\|_{\mathcal{T}_h}^2+\|\nabla v\|_{\mathcal{T}_h}^2+\frac{1}{h_e}|\mu-v|_{\partial\mathcal{T}_h}^2,
\]
and
\[\|(\omega,v,\mu)\|^2_{B,*}:=\|(\omega,v,\mu)\|_B^2+h|\omega\cdot n_e|_{\partial\mathcal{T}_h}^2.
\]
As in \cite{ref8,ref20}, we can read the following stability and boundedness of the bilinear form $B_\mathbf{u}$.
\begin{lemma}[Stability and Boundedness] Assume that $\phi$ and $\phi_h$ are fixed, for all $( \mathbf{u}, p, \lambda) \in  H^k(\mathcal{T}_h) \times L^2(\mathcal{T}_h)\times L^2
(\partial\mathcal{T}_h)$ and $(\omega_h,v_h,  \mu_h) \in \Pi_h \times \Lambda_h \times \Sigma_h$, there  holds
\begin{equation}
\begin{aligned}
&(\textrm{a})\quad| B_{\mathbf{u}}((\mathbf{u},p,\lambda),(\omega_h,v_h,\mu_h))|\leq K^{\ast}\|(\mathbf{u},p,\lambda)\|_{B,\ast}\|(\omega_h,v_h,\mu_h)\|_B,\\
&(\textrm{b})\quad\sup\limits_{(\omega_h,v_h,\mu_h)\in \Pi_h \times \Lambda_h \times \Sigma_h}\frac{B_{\mathbf{u}}((\mathbf{u}_h,p_h,\lambda_h),(\omega_h,v_h,\mu_h))}{\|(\omega_h,v_h,\mu_h)\|_B}\geq K_{\ast}\|(\mathbf{u}_h,p_h,\lambda_h)\|_B,
\end{aligned}
\end{equation}
where $K_{\ast}$ and $K^{\ast}$ denote two positive constants  independent of the mesh size $h$ .
\end{lemma}

Introduce the local $L^2$-projection operators
$\Pi_h$ and $\Pi_e$ as follows:
\begin{equation}
\begin{aligned}
&(p-\Pi_hp,v_h)_K=0,\quad \forall v_h\in P_k(K),\\
&\langle\lambda-\Pi_ep,\mu_h\rangle_e=0,\quad \forall \mu_h\in P_k(e),
\end{aligned}
\end{equation}
 where $K\in\mathcal{T}_h$, $e\in\mathcal{E}_h^i$, $p\in L^2(K)$ and $\lambda\in L^2(e)$.
 
\begin{lemma}[\cite{ref21}]\label{31}
For the local $L^2$-projection operators $\Pi_h$ and $\Pi_e$, there exists the following approximate property
\begin{equation}\label{111}
\begin{aligned}
&\|p-\Pi_hp\|_K\leq Ch^s\|p\|_{s,K},\quad 0\leq s\leq k+1,\\
&\|\nabla(p-\Pi_hp)\|_K\leq Ch^s\|p\|_{s+1,K},\quad 0\leq s\leq k,\\
&\|p-\Pi_hp\|_e+\|p-\Pi_ep\|_e\leq Ch^{s+\frac{1}{2}}\|p\|_{s+1,K},\quad 0\leq s\leq k.
\end{aligned}
\end{equation}
\end{lemma}

The  classical Raviart-Thomas projection operator as in \cite{ref6} is also used 
\begin{equation}
(\mathbf{u}-\Pi^{RT}\mathbf{u},\omega_h)_K=0,\quad\forall\omega_h\in[P_{k-1}(K)]^d,
\end{equation}
and
\begin{equation}
\langle(\mathbf{u}-\Pi^{RT}\mathbf{u})\cdot n_e,\mu_h\rangle_e=0,\quad\forall\mu_h\in P_{k}(e),\quad e\in\partial K.
\end{equation}
We can reach the error estimate as in \cite{ref6}:
\begin{lemma}
For the Raviart-Thomas interpolation $\Pi^{RT}$, the following estimate hold
\begin{equation}
\begin{aligned}
&||\nabla\cdot(\mathbf{u}-\Pi^{RT}\mathbf{u})||_K\leq Ch^s||\nabla\cdot\mathbf{u}||_{s,K},\quad 1\leq s\leq k+1,\\
&||\mathbf{u}-\Pi^{RT}\mathbf{u}||_K+h^{\frac{1}{2}}||u-\Pi^{RT}\mathbf{u}||\leq Ch^s||\mathbf{u}||_{s,K},\quad \frac{1}{2}\leq s\leq k+1.
\end{aligned}
\end{equation}
\end{lemma}
Utlizing the above results element-wise, we can easily get the following  error estimates.
\begin{lemma}\label{10}
If $a(\phi)$ is bounded, there exists the following inequality
\begin{equation}\label{11}
\begin{aligned}
&\|(\mathbf{u}-\Pi^{RT}\mathbf{u},p-\Pi_h p,\lambda-\Pi_e p)\|_{B,\ast}\leq C h^s| p|_{s+1,\mathcal{T}_h},\quad \frac{1}{2}<s\leq k.
\end{aligned}
\end{equation}
\end{lemma}

\begin{remark}
From \eqref{11}, the following estimate holds: for any $1/2\leq s\leq k+1$,
\begin{equation}
\begin{aligned}
&\|\mathbf{u}-\Pi^{RT}\mathbf{u}\|_{K}\leq C h^s\| p\|_{s,K}.
\end{aligned}
\end{equation}
\end{remark}

For the concentration, we introduce another projection operator $\Pi_s$ as follows:
\begin{equation}\label{13}
B_c(c_f-\Pi_s c_f,z_h)+\delta(c_f-\Pi_sc_f,z_h)=0,\quad \forall z_h\in\Psi_h,
\end{equation}
where $\delta$ should be some sufficient large  constant.

As in \cite{ref12}, under the  following inductive hypothesis
\begin{equation}\label{uh}
\| \mathbf{u}_h\|_{L^{\infty}}\leq C_{\mathbf{u}},
\end{equation}
where $C_{\mathbf{u}}$ is a positive constant, we can reach the following estimates:
\begin{equation}\label{err222}
\begin{aligned}
&\| c_f-\Pi_s c_f\|_{s,\mathcal{T}_h}\leq Ch^s\| c_f\|_{s,\mathcal{T}_h},\quad 0\leq s\leq k+1,\\
&\| \frac{\partial(c_f-\Pi_s c_f)}{\partial t}\|_{s,\mathcal{T}_h}\leq Ch^s(| c_f|_{s,\mathcal{T}_h}+| \frac{\partial c_f}{\partial t}|_{s,\mathcal{T}_h}),\quad 0\leq s\leq k+1.
\end{aligned}
\end{equation}
%

The following trace inequalities will be also used to prove the convergence theorem (see Lemma 3.1 in \cite{ref16}).

\begin{lemma}
For $\forall v\in H^1(K)$, the trace inequalities are shown below
\begin{equation}\label{1111}
\begin{aligned}
&\| v\|_{0,e}^2\leq C(h_e^{-1}\| v\|_{0,K}^2+h_e\| v\|_{1,K}^2),\\
&\| \nabla v\cdot n_e\|_{0,e}^2\leq C(h_e^{-1}\| \nabla v\|_{0,K}^2+h_e\|\nabla^2v\|_{0,K}^2).
\end{aligned}
\end{equation}
\end{lemma}

\section{The proof of convergence theorem}

Now, we can complete the proof of our convergence theorem  \ref{thm1}.

\begin{proof}
We firstly give the bound of $\|c_h-c_f\|_{\mathcal{T}_h}$. Set $\xi_c=c_h-\Pi_sc_f$, $\zeta_c=c_f-\Pi_sc_f$. Taking $z_h=\xi_c$ in \eqref{13},  we have
\begin{eqnarray}\label{141}
\begin{aligned}
&\sum\limits_{K\in\mathcal{T}_h}\int_K\frac{(\partial\phi_h\xi_c)}{\partial t}\xi_cdx+\sum\limits_{K\in\mathcal{T}_h}\int_K\phi_h\mathbf{D}(\mathbf{u}_h)\nabla\xi_c\cdot\nabla \xi_cdx
\\
&+\sum\limits_{K\in\mathcal{T}_h}\int_Kf_p{\xi_c}^2dx+J_0^{r}(\xi_c,\xi_c)\\
=&\sum\limits_{K\in\mathcal{T}_h}\int_K\frac{(\partial\phi_h\zeta_c)}{\partial t}\xi_cdx+\sum\limits_{K\in\mathcal{T}_h}\int_K(\phi\mathbf{D}(\mathbf{u})-\phi_h\mathbf{D}(\mathbf{u}_h))\nabla c_f\cdot\nabla\xi_cdx\\
&+2\sum\limits_{e\in\mathcal{E}_h^0}\int_e\{\phi_h\mathbf{D}(\mathbf{u}_h)\nabla\xi_c\cdot n_e\}[\xi_c]ds+\sum\limits_{K\in\mathcal{T}_h}\int_K(\mathbf{u}_h-\mathbf{u})c_f\nabla\xi_cdx\\
&+\sum\limits_{e\in\mathcal{E}_h^i}\int_e\{(\phi_h\mathbf{D}(\mathbf{u}_h)-\phi\mathbf{D}(\mathbf{u}))\nabla c_f)\cdot n_e\}[\xi_c]ds+\sum\limits_{K\in\mathcal{T}_h}\int_K\mathbf{u}\xi_c\nabla\xi_cdx\\
&+\sum\limits_{e\in\mathcal{E}_h^i}\int_e\{(\mathbf{u}-\mathbf{u}_h)c_f\cdot n_e\}[\xi_c]ds-2\sum\limits_{e\in\mathcal{E}_h^i}\int_e\{\mathbf{u}_h\xi_c\cdot n_e\}[\xi_c]ds
\\
&+\sum\limits_{K\in\mathcal{T}_h}\int_K(f_P-\delta)\zeta_c\xi_cdx+ \sum\limits_{K\in\mathcal{T}_h}\int_K\kappa[(1-\phi)c_{f}-(1-\phi_{h})c_{h}]\xi_{c}dx\\
&+ \sum\limits_{K\in\mathcal{T}_h}\int_K\frac{\partial(\phi-\phi_h)}{\partial t}c_f\xi_{c}dx+\sum\limits_{K\in\mathcal{T}_h}\int_K\frac{\partial c_f}{\partial t}(\phi-\phi_h)\xi_{c}dx
\\=&F_1+F_2+\cdots+F_{12}.
\end{aligned}
\end{eqnarray}

Now we estimate the terms on the right hand side of \eqref{141} one by one. Using \eqref{2222} and Lemma \ref{phih}, we can get the following result
\begin{equation}\nonumber
\begin{aligned}
&| F_1|+| F_2|+| F_4|+| F_6|+| F_9|+| F_{10}|+| F_{11}|+| F_{12}|\\
\leq &C\{\|\zeta_c\|_{\mathcal{T}_h}^2+\|\xi_c\|_{\mathcal{T}_h}^2+\|\frac{\partial\zeta_c}{\partial t}\|_{\mathcal{T}_h}^2+\|\frac{\partial(\phi-\phi_h)}{\partial t}\|_{\mathcal{T}_h}^2+\|\phi-\phi_h\|_{\mathcal{T}_h}^2+\|\mathbf{u}-\mathbf{u}_h\|_{\mathcal{T}_h}^2\}+\varepsilon\|\nabla\xi_c\|_{\mathcal{T}_h}^2.
\end{aligned}
\end{equation}

For $F_3$, using \eqref{1111} we have
\begin{equation}\nonumber
\begin{aligned}
| F_3|
\leq&\varepsilon J_0^\gamma(\xi_c,\xi_c)+C\sum\limits_{e\in\mathcal{E}_h^i}\gamma^{-1}h_e\| \nabla\xi_c\|_{L^2(e)}^2
\\\leq&\varepsilon J_0^\gamma(\xi_c,\xi_c)+C_1\gamma^{-1}\| \nabla\xi_c\|_{\mathcal{T}_h}^2.
\end{aligned}
\end{equation}

Next, we estimate $F_5$ with \eqref{2222} and \eqref{1111}
\begin{equation}\nonumber
\begin{aligned}
| F_5|&\leq\|\nabla c_f\|_{L^\infty}\sum\limits_{e\in\mathcal{E}_h^i}\| \mathbf{D}(\mathbf{u})-\mathbf{D}(\mathbf{u}_h)\|_{L^2(e)}\|[\xi_c] \|_{L^(e)}\\
&\leq\varepsilon J_0^\gamma(\xi_c,\xi_c)+C\sum\limits_{e\in\mathcal{E}_h^i}\gamma^{-1}h_e\|\mathbf{u}-\mathbf{u}_h\|_{L^(e)}^2\\
&\leq\varepsilon J_0^\gamma(\xi_c,\xi_c)+C\gamma^{-1}\|\mathbf{u}-\mathbf{u}_h\|_{\mathcal{T}_h}^2.
\end{aligned}
\end{equation}
Using the same technique as above, we can reach
\begin{equation}\nonumber
\begin{aligned}
| F_7|&\leq\varepsilon J_0^\gamma(\xi_c,\xi_c)+C\sum\limits_{e\in\mathcal{E}_h^i}\gamma^{-1}h_e\|\mathbf{u}-\mathbf{u}_h\|_{L^(e)}^2
\\&\leq\varepsilon J_0^\gamma(\xi_c,\xi_c)+C\gamma^{-1}\|\mathbf{u}-\mathbf{u}_h\|_{\mathcal{T}_h}^2,\\
 |F_8|&\leq\varepsilon J_0^\gamma(\xi_c,\xi_c)+\|\xi_c\|_{\mathcal{T}_h}^2,
\end{aligned}
\end{equation}
where $\gamma$ is large enough, $\varepsilon$ is small enough and they satisfy
\begin{equation}\nonumber
\begin{aligned}
&C_1\gamma^{-1}\leq\frac{d_m^{\ast}}{4},\quad\quad\varepsilon\leq \min(\frac{1}{4},\frac{d_m^{\ast}}{4}).
\end{aligned}
\end{equation}

Next, we deal with the first term on the left hand side of \eqref{141}. Since $(\frac{\partial\phi_h}{\partial t}\xi_c,\xi_c)_{\mathcal{T}_h}\geq 0$, we can get
\begin{equation}\nonumber
\begin{aligned}
(\frac{\partial(\phi_h\xi_c)}{\partial t},\xi_c)_{\mathcal{T}_h}\geq\frac{1}{2}\frac{\partial}{\partial t}(\phi_h\xi_c,\xi_c)_{\mathcal{T}_h}.
\end{aligned}
\end{equation}
Substituting these estimates into \eqref{141}, and then integrating it on $t$, we get
\begin{equation}\label{14}
\begin{aligned}
&\|\xi_c\|_{\mathcal{T}_h}^2+\int_{0}^{t}\|\nabla\xi_c\|_{\mathcal{T}_h}^2d\tau+\int_{0}^{t}J_0^\gamma(\xi_c,\xi_c)d\tau\\
\leq&C\int_{0}^{t}(\|\frac{\partial\zeta_c}{\partial t}\|_{\mathcal{T}_h}^2+\| \zeta_c\|_{\mathcal{T}_h}^2+\|\frac{\partial(\phi-\phi_h)}{\partial t}\|_{\mathcal{T}_h}^2+\|\xi_c\|_{\mathcal{T}_h}^2+\|\mathbf{u}-\mathbf{u}_h\|_{\mathcal{T}_h}^2+\|\phi-\phi_h\|_{\mathcal{T}_h}^2)d\tau.
\end{aligned}
\end{equation}

From the above estimate, we need to estimate the bound of $\mathbf{u}-\mathbf{u}_h$. Using the definations of projection operators and $B_{\mathbf{u}}$, we have
\begin{equation}
B_{\mathbf{u}}((\Pi^{RT}\mathbf{u}-\mathbf{u},\Pi_h p-p,\Pi_e p-p),(\omega_h,v_h,\mu_h))
=(a(\phi_h)(\Pi^{RT}\mathbf{u}-\mathbf{u}),\omega_h)_{\mathcal{T}_h}.
\end{equation}
According to the boundedness and stability of  the bilinear form $B_{\mathbf{u}}$,  we have the estimate
\begin{equation}\nonumber
\begin{aligned}
&K_{\ast}\|(\Pi^{RT}\mathbf{u}-\mathbf{u}_h,\Pi_h p-p_h,\Pi_e p-\lambda_h)\|_B\\
\leq&\sup\limits_{(\omega_h,v_h,\mu_h)}\frac{B_{\mathbf{u}}((\Pi^{RT}\mathbf{u}-\mathbf{u}_h,\Pi_h p-p_h,\Pi_e p-\lambda_h),(\omega_h,v_h,\mu_h))}{\| (\omega_h,v_h,\mu_h)\|_B}\\
\leq&C(\| c_f-c_h\|_{\mathcal{T}_h}+\| \Pi^{RT}\mathbf{u}-\mathbf{u}\|_{\mathcal{T}_h}+\| \phi-\phi_h\|_{\mathcal{T}_h}).
\end{aligned}
\end{equation}
Hence we get
\begin{equation}\label{12}
\begin{aligned}
\|(\Pi^{RT}\mathbf{u}-\mathbf{u}_h,\Pi_h p-p_h,\Pi_e p-\lambda_h)\|_B
\leq C\{\| c_f-c_h\|_{\mathcal{T}_h}+\| \Pi^{RT}\mathbf{u}-\mathbf{u}\|_{\mathcal{T}_h}+\| \phi-\phi_h\|_{\mathcal{T}_h}\}.
\end{aligned}
\end{equation}
Using  \eqref{12}, we get the estimate
\begin{equation}\label{322}
\begin{aligned}
&\|(\mathbf{u}-\mathbf{u}_h,\Pi_h p-p_h,\Pi_e p-\lambda_h)\|_B\\
\leq&\|(\Pi^{RT}\mathbf{u}-\mathbf{u}_h,\Pi_h p-p_h,\Pi_e p-\lambda_h)\|_B+\|(\mathbf{u}-\Pi^{RT}\mathbf{u},\Pi_h p-p_h,\Pi_e p-\lambda_h)\|_B\\
\leq&C\{\| c_f-c_h\|_{\mathcal{T}_h}+\|\Pi^{RT}\mathbf{u}-\mathbf{u}\|_{\mathcal{T}_h}+\| \phi-\phi_h\|_{\mathcal{T}_h}\}.
\end{aligned}
\end{equation}

Next, we estimate the boundedness of $\|\phi-\phi_h\|_{\mathcal{T}_h}$. From \eqref{semi}, we can get that
\begin{equation}\label{errphi}
\begin{aligned}
\frac{\partial(\phi-\phi_h)}{\partial t}
&\leq\frac{\alpha \kappa}{\rho_s}[(1-\phi_h)|c_f-c_h|+(\phi_h-\phi)c_f].
\end{aligned}
\end{equation}
So we can get that
\begin{equation}
\begin{aligned}
\|\frac{\partial(\phi-\phi_h)}{\partial t}\|_{\mathcal{T}_h}^2
\leq C(\|\phi-\phi_h\|_{\mathcal{T}_h}^2+\|c_f-c_h\|_{\mathcal{T}_h}^2).
\end{aligned}
\end{equation}
Multiplying \eqref{errphi} by $\phi-\phi_h$ and integrating it over $\Omega$, we will reach that
\begin{equation}\nonumber
\begin{aligned}
\frac{1}{2}\frac{\partial}{\partial t}\|\phi-\phi_h\|_{\mathcal{T}_h}^2\leq C(\|\phi-\phi_h\|_{\mathcal{T}_h}^2+\| c_f-c_h\|_{\mathcal{T}_h}^2).
\end{aligned}
\end{equation}
So we obtain that
\begin{equation}
\begin{aligned}
\|\phi-\phi_h\|_{\mathcal{T}_h}^2\leq C\int_{0}^{t}(\|\xi_c\|_{\mathcal{T}_h}^2+\|\zeta_c\|_{\mathcal{T}_h}^2)d\tau.
\end{aligned}
\end{equation}
Substituting the above estimate into \eqref{14}, and using \eqref{err222}, \eqref{322} and Gronwall's inequality,
we can get the following estimate
\begin{equation}\nonumber
\begin{aligned}
&\| c_f-c_h\|_{\mathcal{T}_h}^2+\|\phi-\phi_h\|_{\mathcal{T}_h}^2+\| \mathbf{u}-\mathbf{u}_h\|_{\mathcal{T}_h}^2\\
\leq&C\int_{0}^{t}(\| \zeta_c\|_{\mathcal{T}_h}^2+\| \frac{\partial\zeta_c}{\partial t}\|_{\mathcal{T}_h}^2+\| \mathbf{u}-\Pi^{RT}\mathbf{u}\|_{\mathcal{T}_h}^2)d\tau\\
\leq&Ch^s(\| c_f\|_{H^1(0,T;H^s(\mathcal{T}_h))}^2+\| p\|_{L^2(0,T;H^{s+1}(\mathcal{T}_h))}^2).
\end{aligned}
\end{equation}
Combined the above estimate with \eqref{322}, we get the second inequality of \eqref{3191}.

It is easily seen that our estimates  are derived under the induction hypothesis \eqref{uh}. Now, we check it. Note that
\begin{equation}\nonumber
\begin{aligned}
\|\mathbf{u}_h\|_{L^{\infty}}&\leq \|\mathbf{u}_h-\Pi^{RT}\mathbf{u}\|_{L^{\infty}}+\|\Pi^{RT}\mathbf{u}-\mathbf{u}\|_{L^{\infty}}+ \|\mathbf{u}\|_{L^{\infty}}\\
&\leq Ch^{s-\frac{d}{2}}+\|\mathbf{u}\|_{L^{\infty}}\leq C_{\mathbf{u}}.\\
\end{aligned}
\end{equation}
Thus, the hypothesis \eqref{uh} holds.

Using the similar technique as in \cite{ref1}, we know that
\begin{equation}\nonumber
\| \Pi_{h}p-p_h\|_{\mathcal{T}_h}\leq C(1+\|\mathbf{u}\|_{L^{\infty}})\| c_h-c_f\|_{\mathcal{T}_h}.
\end{equation}
Using Lemma \ref{31} and  \eqref{3191}(a) , we get  \eqref{3191}(c).
\end{proof}
\section{Numerical Examples}

 In this section, we will test the efficiency of our proposed method by some numerical examples. We firstly use HMFE method for the linear elliptic problem, and then we consider SIPDG method for the convection-diffusion equation. Next, we confirm the convergence rate of our combined method for the  coupled problem. Finally, we apply the combined method to a ``real'' incompressible wormhole problem.

\subsection{Convergence test of HMFE method}
Here we will test the accuracy of the HMFE scheme. the HMFE method is considered for solving the elliptic problem with $RT0-P0$, $RT1-P1$ and $RT2-P2$ elements. The exact solution is taken by $p = \sin\pi x \sin\pi y$ in $[0,1]\times[0,1]$ and  $\mathbf{u} = -\nabla p$,
respectively. For different mesh size $h=1/8, 1/16, 1/32,1/64$, a convergence study is presented.   The $L^2$-norm
errors and convergence accuracies are shown in Tables \ref{tab1}-\ref{tab21}. As seen in these tables, the optimal convergence
rates for pressure and velocity are evaluated.
{\centering
\begin{table}[ht]
\caption{Numerical results for $p$ and $\mathbf{u}$ with $RT0-P0$.}\label{tab1}
\begin{tabular}{p{1cm}|p{2.5cm}p{1.5cm}|p{2.5cm}p{1.5cm}}
  \hline
  $h$ & $\|p-p_{h}\|_{L^{2}}$  & rates & $\|\mathbf{u}-\mathbf{u}_{h}\|_{L^{2}}$   & rates\\
  \hline
  1/8 & 7.1830e-02 & *  & 2.4473e-02 & * \\
  1/16 &3.5977e-02 & 0.9975  & 1.2508e-02 & 0.9684 \\
  1/32 &1.7992e-02 & 0.9997  & 6.2691e-03 & 0.9965\\
  1/64 & 8.9969e-03 & 0.9985  & 3.1335e-04 & 1.0005\\
  \hline
\end{tabular}
\end{table}
}

{\centering
\begin{table}[ht]
\caption{Numerical results for $p$ and $\mathbf{u}$ with $RT1-P1$.}\label{tab2}
\begin{tabular}{p{1cm}|p{2.5cm}p{1.5cm}|p{2.5cm}p{1.5cm}}
  \hline
  $h$ & $\|p-p_{h}\|_{L^{2}}$  & rates & $\|\mathbf{u}-\mathbf{u}_{h}\|_{L^{2}}$   & rates\\
  \hline
  1/8 & 2.7875e-02 & * & 1.1113e-02 & * \\
  1/16 & 7.1654e-03& 1.9599 & 2.8341e-03 & 1.9713 \\
  1/32 & 1.8225e-03 & 1.9752 & 7.1763e-04 & 1.9816 \\
  1/64 & 4.6070e-04 & 1.9840 & 1.8099e-04 & 1.9874 \\  \hline
\end{tabular}
\end{table}}

{\centering
\begin{table}[ht]
\caption{Numerical results for $p$ and $\mathbf{u}$ with $RT2-P2$.}\label{tab21}
\begin{tabular}{p{1cm}|p{2.5cm}p{1.5cm}|p{2.5cm}p{1.5cm}}
  \hline
  $h$ & $\|p-p_{h}\|_{L^{2}}$  & rates & $\|\mathbf{u}-\mathbf{u}_{h}\|_{L^{2}}$   & rates\\
  \hline
  1/8 &4.4473e-04 &*            &2.5267e-02&*\\
  1/16 &5.8189e-05 &2.9341 &2.9358e-03&3.1054\\
  1/32 &7.3831e-06 &2.9785&3.3422e-04&3.1349\\
  1/64 &9.2786e-07 &2.9922 &3.9131e-05&3.0944\\  
  \hline
\end{tabular}
\end{table}}

\subsection{Convergence test of SIPDG method}
Here we first test the convergent accuracy of the SIPDG method for convection-diffusion equation
\[
\f{\p c_{f}}{\p t}+\nabla\cdot(\mathbf{u}c_{f}-D\nabla c_{f})=f
\]
with homogeneous and nonhomogeneous boundary value conditions. For this purpose, we take the two different exact solutions respectively as
\[
c_{f}=e^{-t}\sin\pi x\sin\pi y, \quad\textrm{and}
\quad
c_{f}=e^{- y^2  - x- t}, \quad (x,y)\in[0,1]\times[0,1].
\]
The velocity function $\mathbf{u}=[-y,x]$ and the diffusion coefficient $D=1.0$. The initial-boundary  conditions and the right hand side term can be computed by the exact solutions. For the practical computation, the first-order Euler backward difference scheme in time is used and $L^{2}$-projection of  the initial condition is also used. Setting $T=1.0$ and time size $\Delta t=1e-3$, for different mesh size, we give some numerical results with $P1$  discontinuous finite element space in Table \ref{DG1}.  
These numerical results show that SIPDG method has the optimal convergence rates in $L^{2}$-norm for both homogeneous and nonhomogeneous boundary conditions.

{\centering
\begin{table}[ht]
\caption{Numerical results with $P1$ element for homogeneous and nonhomogeneous boundary cases.}\label{DG1}
\begin{tabular}{p{1cm}|p{2.5cm}p{1.5cm}|p{2.5cm}p{1.5cm}}
  \hline
  h&\multicolumn{2}{|c}{homogenous}&\multicolumn{2}{|c}{nonhomogenous}\\
  \cline{2-5}
  & $L^{2}$ error & rates & $L^{2}$ error & rates\\
  \hline
  1/8   &1.1599e-00 &* &2.2602e-02 &*\\
  1/16 &3.0175e-01 &1.9425&5.8884e-03 &1.9405 \\
  1/32 &7.6512e-02 &1.9796&1.5668e-03 &1.9100\\
 1/ 64 &1.9192e-02&1.9952&4.1542e-04&1.9152 \\
  \hline
\end{tabular}
\end{table}}

In addition, we also consider our method for the porosity and the concentration.  Initial-boundary conditions can be given by the exact solutions
\begin{equation}
\begin{aligned}
&c_f(x,y,t)=\f{2\epsilon^2}{2\epsilon^2+4Dt}\exp\{-\frac{(x\cos4t+y\sin4t+0.2)^2+(-x\sin4t+y\cos4t)^2}{2\epsilon^2+4Dt}\},\\
&\phi(x,y,t)=0.5+0.4\sin(x+t)\sin(y+t),\quad \textrm{in}\quad \Omega= [0, 1]\times[0, 1].
\end{aligned}
\end{equation}
The other parameters are taken as:
\begin{equation}
c_I=k_c=k_s=a_0=\frac{\alpha}{\rho_s}=1,\quad D=0.1,\quad \epsilon=0.1.
\end{equation}
Here we still use the first-order backward Euler scheme in time,  and take time step $\Delta t = 1e-3$.
The computational results at $T = 1.0$ are shown as in Tables \ref{tab3} and \ref{tab4} with the uniform triangular meshes  $h = 1/8, 1/16, 1/32, 1/64, 1/128$.
From these tables, we can get the optimal convergence rates in $L^2$-norm with $P1$ and $P2$ discontinuous elements.

{\centering
\begin{table}[ht]
\caption{Numerical results with $P1$ element for $c_f$ and $\phi$.}\label{tab3}
\begin{tabular}{p{1cm}|p{2.5cm}p{1.5cm}|p{2.5cm}p{1.5cm}}
  \hline
  $h$ & $\|c_f-c_{h}\|_{L^{2}}$  & rates & $\|\phi-\phi_{h}\|_{L^{2}}$   & rates\\
  \hline
  1/8 & 7.0710e-02 & * & 8.2931e-02 & * \\
  1/16 & 1.7289e-02 & 2.0320 & 2.4605e-02 & 1.7529 \\
  1/32 & 4.2388e-03 & 2.0281 & 6.7051e-03 & 1.8756\\
  1/64 & 1.0472e-03 & 2.0171 & 1.7132e-03 & 1.9686\\
  1/128 & 2.6026e-04 & 2.0085 & 4.3005e-04 & 1.9941 \\
  \hline
  \end{tabular}
\end{table}}

{\centering
\begin{table}[ht]
\caption{Numerical results  with $P2$ element for $c_f$ and $\phi$.}\label{tab4}
\begin{tabular}{p{1cm}|p{2.5cm}p{1.5cm}|p{2.5cm}p{1.5cm}}
  \hline
  $h$ & $\|c_f-c_{h}\|_{L^{2}}$  &  rates & $\|\phi-\phi_{h}\|_{L^{2}}$   & rates\\
  \hline
  1/8 & 5.2783e-03 & * & 2.7432e-02 & * \\
  1/16 & 6.3176e-04 & 3.0626 & 4.0521e-03 & 2.7591 \\
  1/32 & 7.5796e-05 & 3.0592 & 5.3296e-04 & 2.9266\\
  1/64 & 9.1162e-06 & 3.6314 & 6.7731e-05 & 2.9761\\
  1/128 & 1.0039e-06 & 3.1828 & 8.0221e-06 & 3.0429 \\
  \hline
\end{tabular}
\end{table}}

\subsection{Convergence test of the combined method}
In this experiment, we will  show the convergence of our combined method. Here the analytic solution  in $\Omega=[0,1]\times[0,1]$ is given  as in \cite{rui2}
\begin{equation}\nonumber
\begin{aligned}
&p(x,y,t)=t\cos\pi x\cos\pi y,\\
&c_f(x,y,t)=tx^2(1-x)^2y^2(1-y)^2,\\
&\phi(x,y,t)=1-e^{-\frac{1}{80}t^2x^2(1-x)^2y^2(1-y)^2e^{x+y+1}-(x+y+1)}.
\end{aligned}
\end{equation}
The parameters are taken as
\begin{equation}\nonumber
\begin{aligned}
&\mathbf{D}=10^{-2}\mathbf{I},\quad k_0=1,\quad a_0=0.5,\quad\rho_s=10,\quad\alpha=1,\quad k_c=k_s=1,\quad\mu=f_I=1,
\end{aligned}
\end{equation}
where I is an identity matrix. And choosing $T=1.0$ and time step $\Delta t=h^2$, we give some numerical results with $RT1-P1$ element and $P1$ discontinuous element in Tables \ref{FEDG1} and \ref{FEDG2}. We can easily find that our combined method is of second-order accuracy in $L^{2}$-norm, which is coincided with our theoretical analysis.

{\centering\begin{table}[ht]
\caption{Numerical results  for $c_f$ and $\phi$ with $P1$ element.}\label{FEDG1}
\begin{tabular}{p{1cm}|p{2.5cm}p{1.5cm}|p{2.5cm}p{1.5cm}}
  \hline
  $h$ & $\|c_f-c_{h}\|_{L^{2}}$  & rates & $\|\phi-\phi_{h}\|_{L^{2}}$   & rates\\
  \hline
  1/8 & 1.1109e-03 & * & 2.4572e-02 & * \\
  1/16 & 2.9657e-04 & 1.9053 & 6.8963e-03 & 1.8331 \\
  1/32 & 7.6954e-05 & 1.9463 & 1.7291e-03 & 1.9958\\
  1/64 & 1.9001e-05 & 2.0181 & 4.2112e-04 & 2.0377\\
  1/128 & 4.4123e-06 & 2.1065 & 1.0021e-04 & 2.0712 \\
  \hline
\end{tabular}
\end{table}}
{\centering
\begin{table}[ht]
\caption{Numerical results for $\mathbf{u}$ and $p$  with $RT1-P1$ element.}\label{FEDG2}
\begin{tabular}{p{1cm}|p{2.5cm}p{1.5cm}|p{2.5cm}p{1.5cm}}
  \hline
  $h$ & $\|\mathbf{u}-\mathbf{u}_{h}\|_{L^{2}}$  & rates & $\|p-p_{h}\|_{L^{2}}$   & rates\\
  \hline
  1/8 & 2.4932e-03 & * & 6.2173e-03 & * \\
  1/16 & 6.2776e-04 & 1.9897 & 1.7321e-03 & 1.8438 \\
  1/32 & 1.7290e-04 & 1.8603 & 3.9021e-04 & 2.1502\\
  1/64 & 4.2003e-05 & 2.0414 & 9.7001e-05 & 2.0082\\
  1/128 & 1.0010e-05 & 2.0691 & 2.5231e-05 & 1.9428 \\
  \hline
\end{tabular}
\end{table}}

\subsection{Simulation for a ``real'' incompressible wormhole propagation}
In this experiment, a $0.2$-meter computational domain is considered,  and the first-order Euler backward time discretization is used. We set a singular area on the middle of the left boundary with space size to be $0.01$-meter and time size to be $1e-4$ to observe the phenomenon of wormhole propagation. The initial values and the parameters in the porous medium are taken as in Table \ref{tab5}. Initial concentration of acid and initial porosity of rock in this domain are set to be $c_0 = 0$ and $\phi_0 = 0.2$, respectively. The top and bottom boundaries of the domain are impermeable.
\newpage

{\centering
\begin{table}[ht]
\caption{The properties of acid flow and porous medium.}\label{tab5}
\begin{tabular}{l| l}
\hline
Properties & Value\\
\hline
 the viscosity of fluid $(\mu)$ & $1$$Pa\cdot s$\\
 the injection flow rate $(f_I)$ & 4.5\\
 the production flow rate $(f_P)$ & 2.5\\
 the dispersion tensor $(\mathbf{D})$ &0.01\\
 the local mass-transfer coefficient $(k_c)$ &1$m/s$\\
 the density of the rock $(\rho_s)$ & 2000$kg/m^2$\\
 the dissolving constant of the acid $(\alpha)$ &0.1$kg/mole$\\
 the kinetic constant for reaction $(k_s)$ &10$m/s$\\
 the initial interfacial area available for reaction $(a_0)$ &0.2$m^{-1}$\\
\hline
\end{tabular}
\end{table}
}

\begin{figure}[ht]
\centering
\subfigure[$\phi$ at $T=10$]{
\includegraphics [width=3cm,height=3cm]{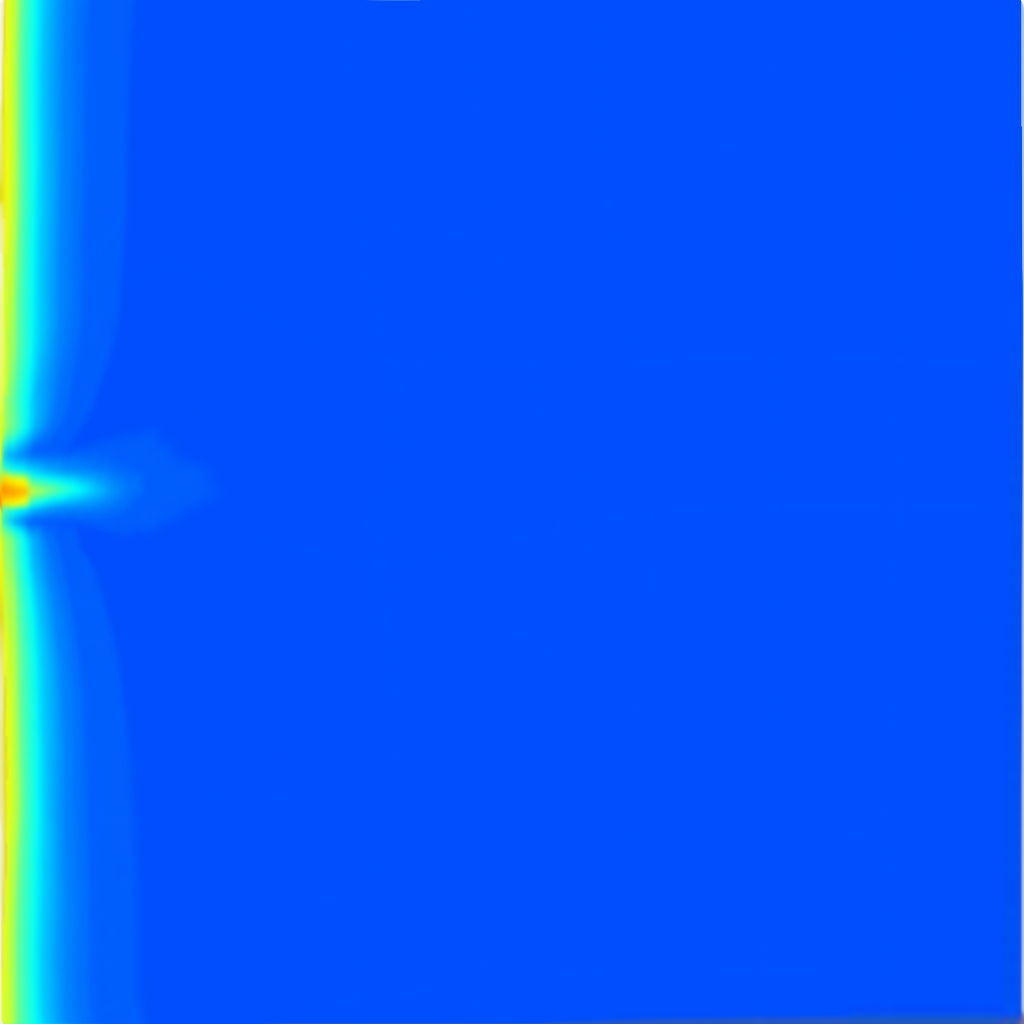}}
\subfigure[$\phi$ at $T=20$]{
\includegraphics [width=3cm,height=3cm]{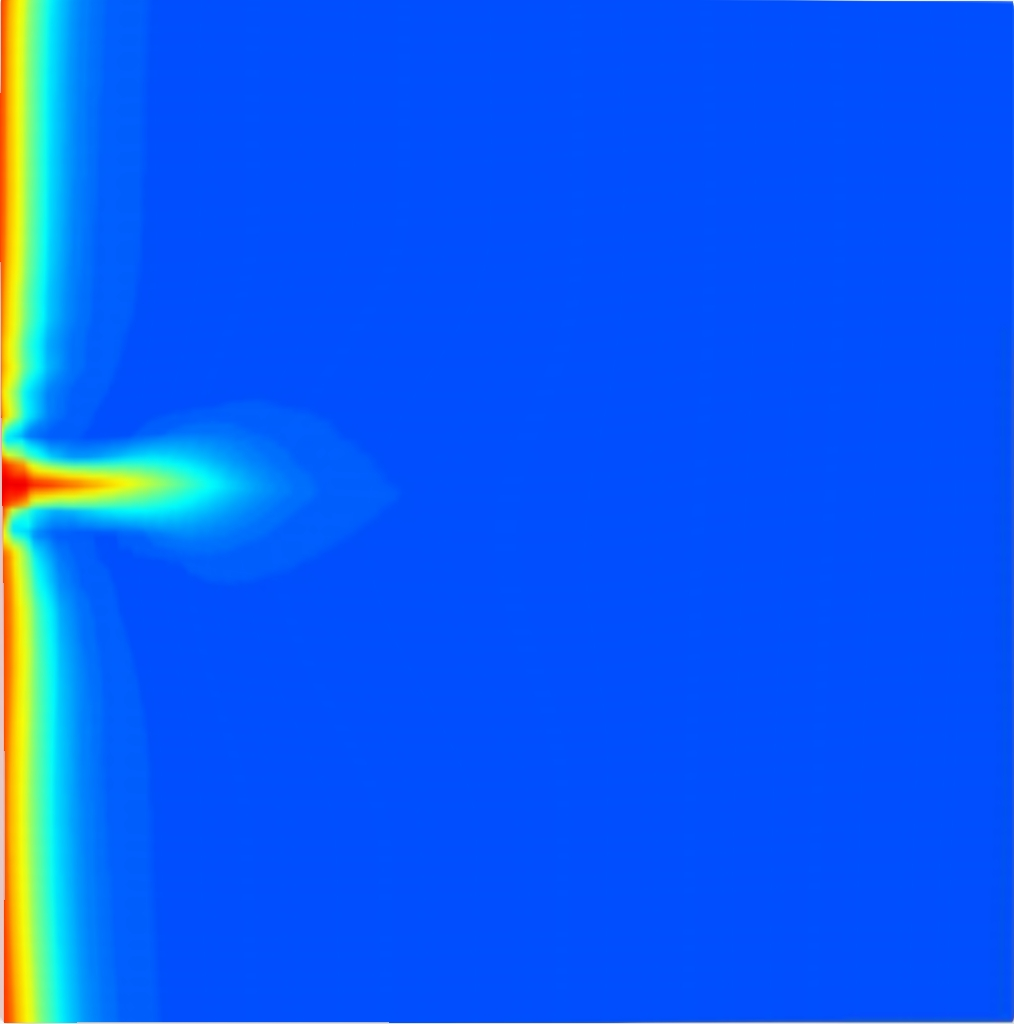}}
\subfigure[$\phi$ at $T=30$]{
\includegraphics [width=3cm,height=3cm]{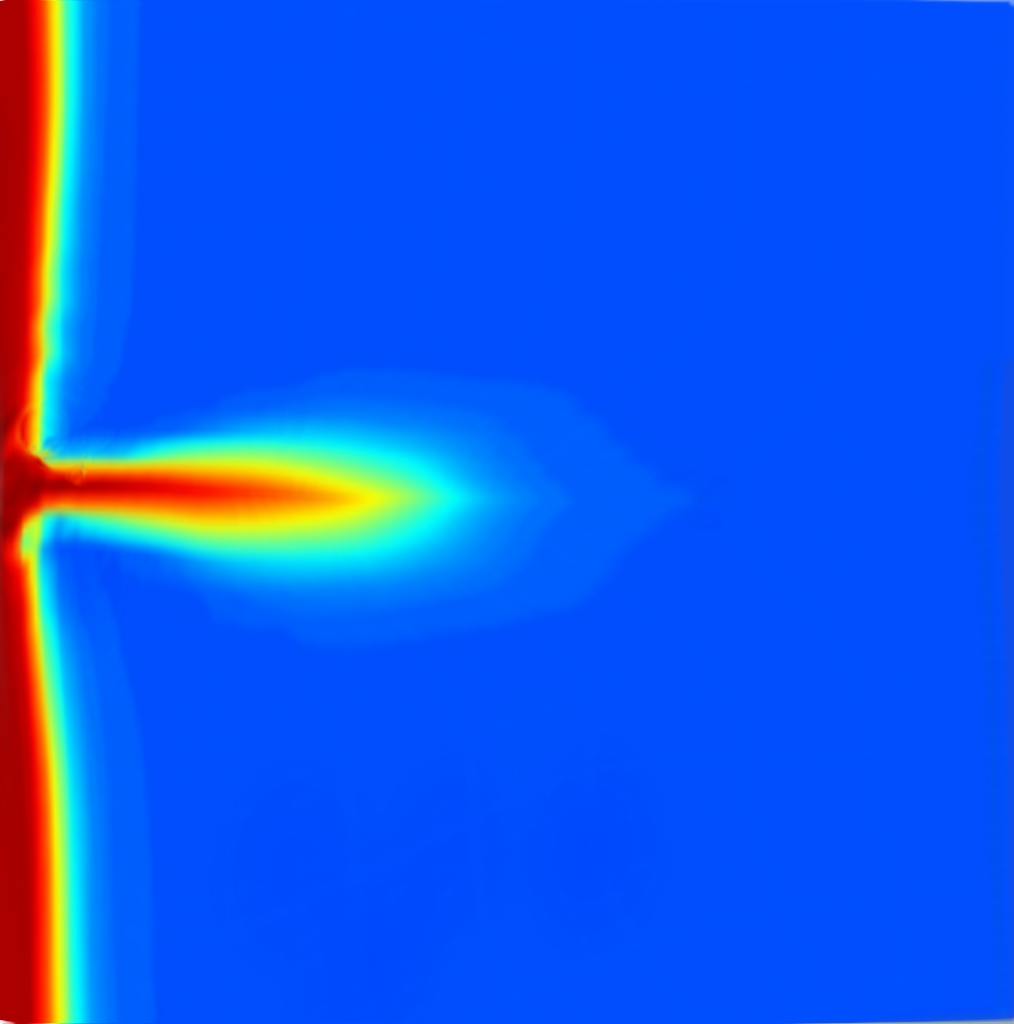}}
\subfigure{
\includegraphics [width=0.5cm,height=3cm]{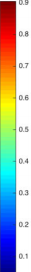}}
\caption{Porosity of rock at the different time steps.}\label{fig1}
\end{figure}

\begin{figure}[ht]
\centering
\subfigure[$c_f$ at $T=10$]{
\includegraphics [width=3cm,height=3cm]{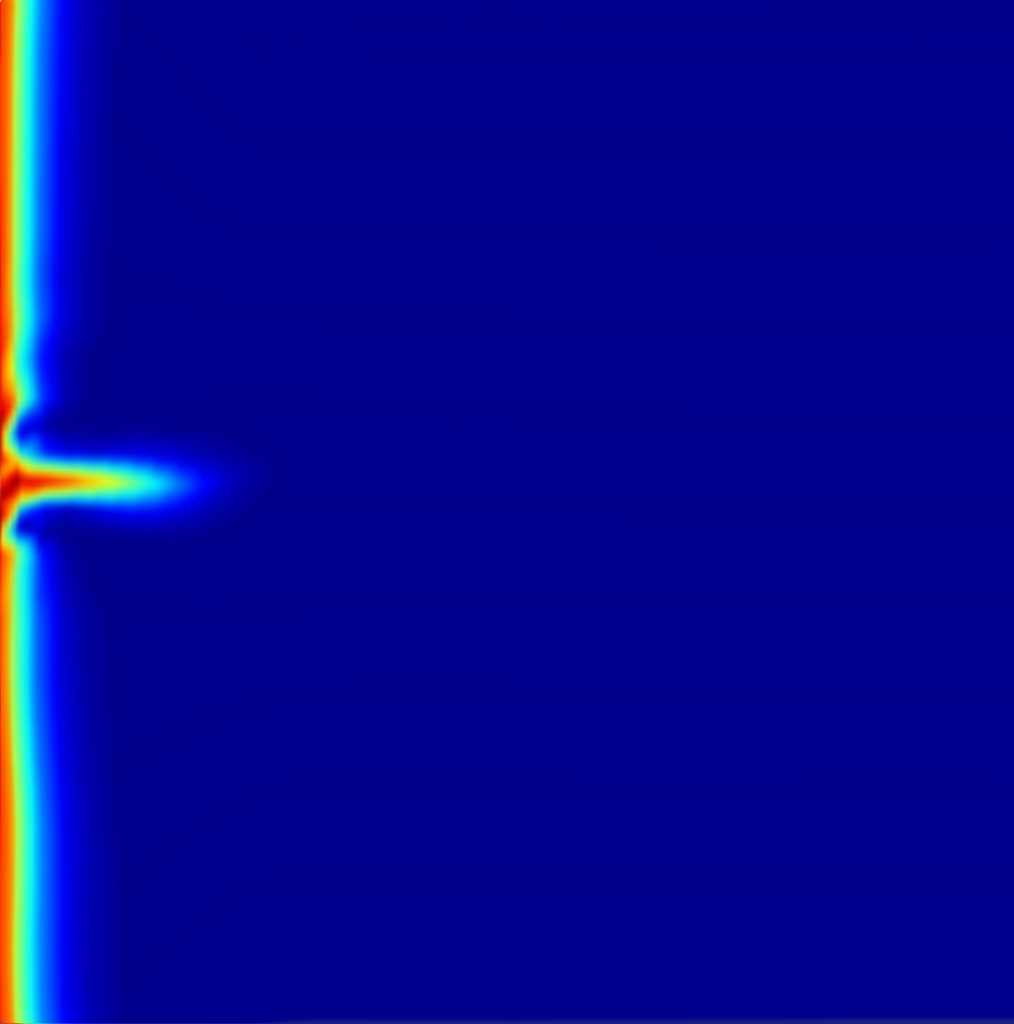}}
\subfigure[$c_f$ at $T=20$]{
\includegraphics [width=3cm,height=3cm]{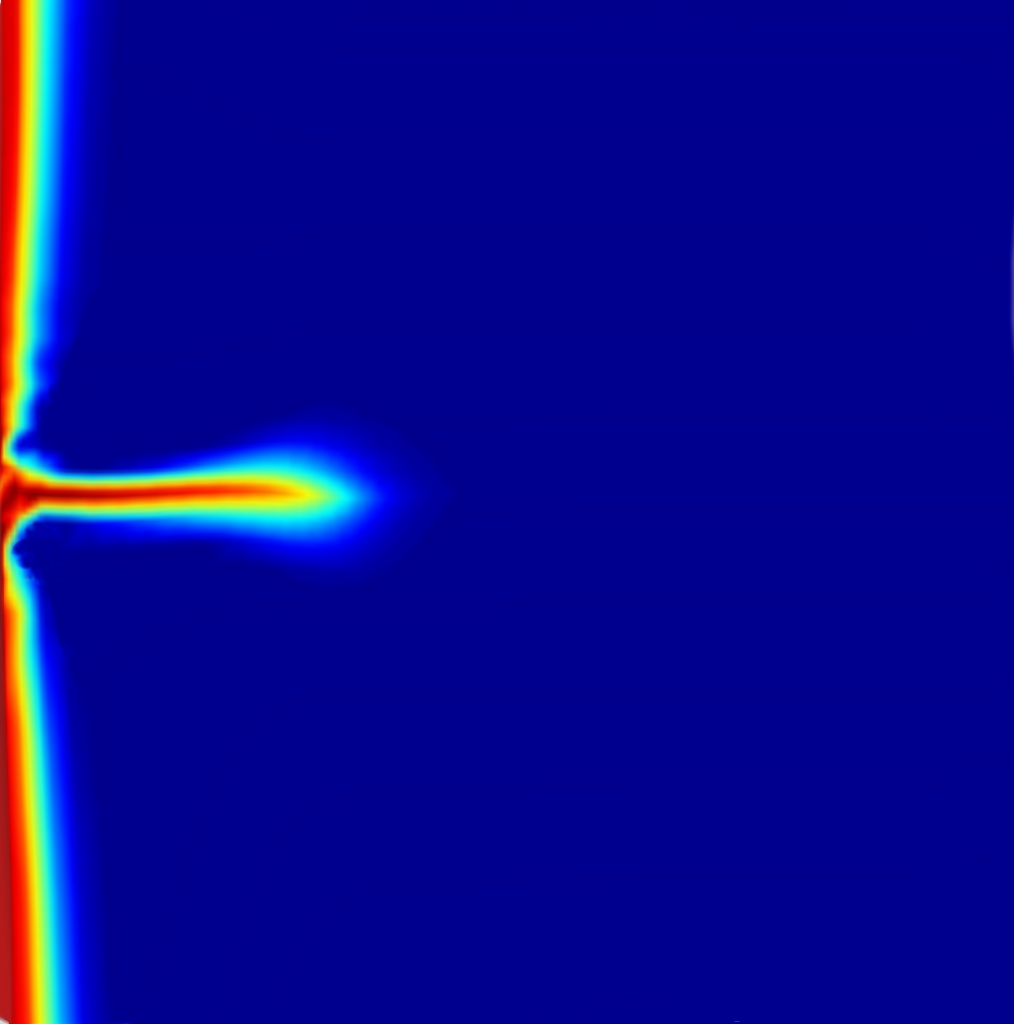}}
\subfigure[$c_f$ at $T=30$]{
\includegraphics [width=3cm,height=3cm]{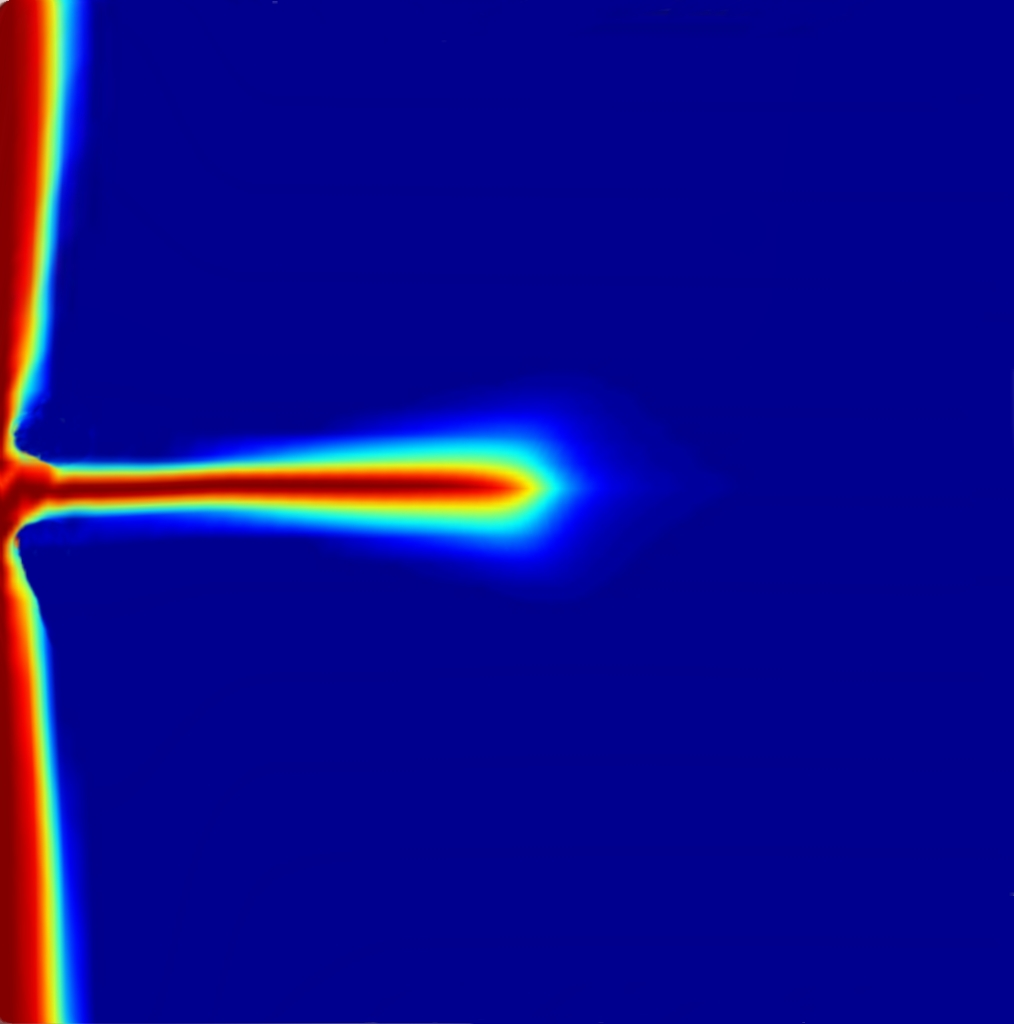}}
\subfigure{
\includegraphics [width=0.5cm,height=3cm]{ex3.png}}
\caption{Concentration of acid at the different time steps.}\label{fig2}
\end{figure}

The numerical results of the concentration and porosity at different time are shown in
Figures \ref{fig1} and \ref{fig2}. From these figures, we can observe $c_f , \phi \in [0, 1]$ and the phenomenon of wormhole propagation, which shows the effectiveness of the combined  method.

%
%
%

%
%

\end{document}